\documentclass[1pt]{article}
\usepackage{amsmath, amssymb, amsthm, amsfonts, cases}
\usepackage{mathrsfs}
\usepackage{url}
\usepackage{authblk}
\usepackage[usenames]{color}
\usepackage{geometry}

\newtheorem{thm}{Theorem}[section]

\newtheorem{ass}{Assumption}[section]

\renewcommand{\theequation}{\arabic{equation}}

\allowdisplaybreaks

 \makeatother
\begin{document}

\title{  Partial Information Stochastic Differential Games  for Backward Stochastic
 Systems Driven By L\'{e}vy Processes 
 \thanks{This work was supported by the Natural Science Foundation of Zhejiang Province
for Distinguished Young Scholar  (No.LR15A010001),  and the National Natural
Science Foundation of China (No.11471079, 11301177) }}

\date{}

   \author{ Fu  Zhang$^{a}$ \hspace{1cm}
   Qingxin Meng$^{b}$ \thanks{Corresponding author.   E-mail: mqx@zjhu.edu.cn}
 \hspace{1cm}  Maoning Tang$^{b}$
\hspace{1cm}
\\ \small{$^{a}$Department of Mathematical Sciences,  University of Shanghai for Science and Technology, Shanghai 200433, China}
\\\small{$^{b}$Department of Mathematics, Huzhou University, Zhejiang 313000, China}}

\maketitle
\begin{abstract}
In this paper, we consider  a partial information  two-person zero-sum stochastic
differential game problem where the system is governed by  a  backward
stochastic differential  equation   driven by Teugels martingales associated with a
L\'{e}vy process and an independent Brownian motion.
One sufficient (a verification
theorem) and one necessary conditions for the existence of optimal controls are proved. To
illustrate the general results, a linear quadratic stochastic differential game problem  is discussed.
\end{abstract}


\maketitle

\section{INTRODUCTION}
\setcounter{equation}{0}
\renewcommand{\theequation}{1.\arabic{equation}}

In this paper,  we consider  a partial information  two-person zero-sum stochastic
differential game problem where the system is governed by the following
nonlinear backward stochastic differential equation (BSDE)
{\small \begin{equation} \label{eq:2.1}
\begin{array}{ll}
&y(t)=\xi+\displaystyle\int_t^Tf(s,y{(s)},q(s),z(s),u_1(t),u_2(t))ds\\&
-\displaystyle\sum_{i=1}^d\int_t^Tq^i(s)dW^i(s)
-\displaystyle\sum_{i=1}^\infty\int_t^T z^i(s)dH^i(s),\ t\in[0,T]
\end{array}
\end {equation}
with the cost functional
\begin{equation}\label{eq:2.2}
J(u(\cdot))=E\displaystyle\bigg[\int_0^Tl(t,y(t),q(t),z(t),u_1(t),u_2(t))dt+\phi(y(0))\bigg],
\end{equation}}

where $\{W(t),0\leq t \leq T\}$ is a standard $d$-dimensional
Brownian motion and  Teugel's martingales $H(t)=\{H^i(t)\}_{i=1}^\infty ,0\leq t \leq T$
associated with a  L\'{e}vy processes $\{L(t), 0\leq t \leq T\}$ (see section 2 for more details). The filtration generated by the
underlying Brownian motion $W$ and  the
L\'{e}vy process $\{L(t), 0\leq t \leq T\}$ is denoted by $\{{\cal F}_t\}_{0\leq t\leq T}$ .

In the above, the processes $u_1(\cdot)$ and
$u_2(\cdot)$ are our  open-loop
control processes which present the controls of the two players.
Let ${
U}_1\subset R^{k_1}$ and ${U}_2\subset R^{k_2}$  are two given nonempty convex sets.
Under many situations under which the
full information ${\cal F}_t$ is inaccessible for players, ones can only observe a partial information.
 For
this,
an admissible control process $u_i(\cdot)$ for the player $i$  is defined
as a ${{\cal G}}_t$-predictable process with values in $U_i$ s.t.
$E\displaystyle\int_0^T|u_i(t)|^2dt<+\infty$ where $i=1,2$.  Here
${\cal G}_t\subseteq {\cal F}_t $ for all $t\in [0,~T]$ is a
given subfiltration representing the information available to the
controller at time t. For ~ example, we could have
$\varepsilon_t={\cal F}_{(t-\delta)^+},~~~t\in [0,~T],$
 where $\delta>0$ is a fixed delay of information.

 The set of all
 admissible  open-loop controls $u_i(\cdot)$ for the player $i$ is
 denoted by ${\cal A}_i,(i=1,2).$    ${\cal A}_1\times {\cal A}_2$
is called  the set of open-loop admissible controls for the players.
We denote the strong solution of (\ref{eq:2.1})  by
$(y^{u_1,u_2}(\cdot), q^{u_1,u_2}(\cdot), z^{u_1,u_2}(\cdot))$, or $(y(\cdot), q(\cdot),
z(\cdot))$ if its dependence on admissible control $({u}_1(\cdot), {u}_2(\cdot))$ is
clear from context. Then we call $(y(\cdot), q(\cdot), z(\cdot))$
the state processes corresponding to the control process $(\bar{u}_1(\cdot), \bar{u}_2(\cdot))$
and call $({u}_1(\cdot),{u}_2(\cdot), y(\cdot), q(\cdot), z(\cdot))$ the admissible
pair.

The partial information  two-person zero-sum stochastic differential
game amounts to determining an admissible open-loop control
$(\bar{u}_1(\cdot), \bar{u}_2(\cdot))$ over  ${\cal A}_1\times{\cal A}_2$
such that
\begin{equation}\label{c1}
J(\bar{u}_1(\cdot), \bar{u}_2(\cdot))=\sup_{u_2(\cdot)\in {\cal
A}_2}\big(\inf_{u_1(\cdot)\in {\cal
A}_1}J(u_1(\cdot),u_2(\cdot))\big).
\end{equation}
Such control $(\bar{u}_1(\cdot), \bar{u}_2(\cdot))$ is  called an
optimal control. The corresponding strong solution  $(\bar{y}(\cdot),\bar{q}(\cdot),\bar{z}(\cdot)) $ of (1.1)
is called   the optimal state process.  Then
$(\bar{u}_1(\cdot),\bar{u}_2(\cdot); \bar{y}(\cdot), \bar{q}(\cdot), \bar{z}(\cdot))$ is
called an optimal pair.  For this partial stochastic differential
game problem, we denote by Problem (P).

Roughly speaking, for the zero-sum  differential game, Player I seek
control $\bar {u}_1(\cdot)$ to minimize (1.2), and Player II seek
control $\bar{u}_2(\cdot)$ to maximize (1.2). Therefore, (1.2)
represents the cost for Player I and the payoff for Player II. There
are basically two types of controls for both players: open-loop
controls and closed- loop controls. In this paper, we concentrate on
the open-loop case. If $(\bar{u}_1(\cdot),\bar{u}_2(\cdot))$ is an
optimal open-loop control, then
\begin{eqnarray}
J(\bar{u}_1(\cdot), {u}_2(\cdot))\leq J(\bar{u}_1(\cdot),
\bar{u}_2(\cdot))\leq J(u_1(\cdot), \bar{u}_2(\cdot))
\end{eqnarray}
for all admissible open-loop controls
$(u_1(\cdot),{u}_2(\cdot))\in {\cal A}_1\times{\cal A}_2$. We
refer to $(\bar{u}_1(\cdot), \bar{u}_2(\cdot))$ as an open-loop
saddle point of Problem (P).

Game theory had been an active area of research and a useful tool in
many applications, particularly in biology and economic. For the partial information two-person zero-sum stochastic differential
games, the objective is to find an saddle point for which the
controller has less information than the complete information
filtration $\{{\cal F}_t\}_{t\geq 0}$.  Recently, An and
${\O}$ksendal [1] established a maximum principle for stochastic
differential games of forward systems with Poisson jumps for the type of partial
information in our paper. Moreover,  We refer to [2-3] and the references therein for more associated
results  on the partial information  stochastic differential games.

In 2000, Nualart and Schoutens[4]
 got a martingale representation theorem for a type of
L\'{e}vy processes through Teugel's martingales, where Teugel's
martingales are a family of pairwise strongly orthonormal
martingales associated with L\'{e}vy processes. Later, they [5]
proved the existence and uniqueness theory of BSDE driven by
Teugel's martingales. The above results are further extended to the
one-dimensional BSDE driven by Teugel's martingales and an
independent multi-dimensional Brownian motion by Bahlali et al[6].

Since the theory of BSDE driven by Teugel¡¯s
martingales and an independent Brownian motion
is established, it is natural to apply the theory to
the stochastic optimal control problem. Now  the full information
 stochastic optimal
 control problem related to
Teugel's martingales has  been studied by many authors. For example,
the stochastic linear-quadratic problem with L\'{e}vy processes was
studied by Mitsui and Tabata[7]. Motivated
by[7], Meng and Tang[8] studied the general full information stochastic optimal
control problem for the forward stochastic systems driven by
Teugel's martingales and an independent multi-dimensional Brownian
motion and prove the corresponding stochastic maximum principle.
Furthermore, Tang and Zhang [9] extend [8] to   the
Backward stochastic systems  and obtain the the corresponding
stochastic maximum principle. For the case of the partial information, in 2012,
Bahlali, Khelfallah  and  Mezerdi [10]
studied the stochastic control problem for forward
system and obtain the corresponding stochastic maximum principle. In the meantime, Meng, Zhang  and Tang
[11] extend  [9] to the the partial information stochastic optimal control problem of backward stochastic systems
and get the corresponding optimality conditions.

But to our best knowledge,  there is
no discussion on the
 partial information  stochastic differential games
 for the system driven by Teugel martingales and an independent Brownian motion, which motives us to write
this paper.
The main purpose of this paper  is to establish partial information
necessary and sufficient conditions for optimality for Problem (P) by  using of  the results in [9]. The
results obtained in this paper can be considered as a generalization
of stochastic optimal control problem  to the two-person zero-sum
case.  As an application, a  two-person zero-sum
stochastic differential game
of linear backward stochastic differential equations with a quadratic
cost criteria under partial information is discussed and the
optimal control is characterized explicitly by adjoint
processes.

The rest of this paper is organized as follows. We
introduce useful notation  and give needed assumptions in section 2. Section 3 is
devoted to present  sufficient conditions for the existence of the optimal control control.
In section 4,  we establish  necessary conditions of optimality.
In the last section, a linear quadratic stochastic differential game problem
 is  solved by applying the theoretical results.

\section{Preliminaries and assumptions}
 \setcounter{equation}{0}
\renewcommand{\theequation}{2.\arabic{equation}}
Let $(\Omega, {F},\{{\cal F}_t\}_{0\leq t\leq T}, P)$ be a complete
probability space. The filtration $\{{\cal F}_t\}_{0\leq t\leq T}$
is right-continuous and generated by a $d$-dimensional standard
Brownian motion $\{W(t), 0\leq t\leq T\}$ and a one-dimensional
L\'{e}vy process $\{L(t), 0\leq t \leq T\}$. It is known that $L(t)$
has a characteristic function of the form: {\small $Ee^{i\theta
L(t)}=\exp\bigg[ia\theta
t-{1\over2}\sigma^2\theta^2t+t\int_{{R}^1}(e^{i\theta x}-1-i\theta x
I_{\{|x|<1\}})v(dx)\bigg],$} where $a\in {R}^1$, $\sigma>0$ and $v$
is a measure on ${R}^1$ satisfying (i)$\displaystyle
\int_0^T(1\wedge x^2)v(dx)<\infty$ and (ii) there exists
$\varepsilon >0$ and $\lambda >0$, s.t. $\displaystyle
\int_{\{-\varepsilon, \varepsilon\}^c} e^{\lambda |x|}v(dx)<\infty$.
These settings imply that the random variables $L(t)$ have moments
of all orders. We denote by $\{H^i(t), 0\leq t \leq
T\}_{i=1}^\infty$ the Teugel's martingales
associated with the L\'{e}vy process $\{L(t),0\leq t \leq T\}$. 
Here $H^i(t)$ is given by
$$
H^i(t)=c_{i,i}Y^{(i)}(t)+c_{i,i-1}Y^{(i-1)}(t)+\cdots+c_{i,1}Y^{(1)}(t),
$$
where $Y^{(i)}(t)=L^{(i)}(t)-E[L^{(i)}(t)]$ for all $i\geq 1$,
$L^{(i)}(t)$ are so called power-jump processes with
$L^{(1)}(t)=L(t)$, $L^{(i)}(t)=\displaystyle\sum_{0<s\leq t}(\Delta
L(s))^i$ for $i\geq 2$ and the coefficients $c_{i,j}$ correspond to
the orthonormalization of polynomials $1,x, x^2,\cdots$ w.r.t. the
measure $\mu(dx)=x^2v(dx)+\sigma^2\delta_0(dx)$. The Teugel's
martingales $\{H^i(t)\}_{i=1}^\infty$ are pathwise strongly
orthogonal 
and their predictable quadratic variation processes are given by
$$\langle H^{(i)}(t), H^{(j)}(t)\rangle=\delta_{ij}t.$$
For more details of Teugel's martingales, we invite the reader to
consult Nualart and Schoutens [10,11]. Denote by $\cal{P}$  the
predictable sub-$\sigma$ field of ${\cal B}([0, T])\times \cal{F}$,
then we introduce the following notation used throughout this paper.

In the following, we introduce some basic spaces.

$\bullet$~~$H$: a Hilbert space with norm $\|\cdot\|_H$.

$\bullet$~~$\langle\alpha,\beta\rangle:$ the inner product in $
{R}^n, \forall \alpha,\beta\in {R}^n.$

$\bullet$~~$|\alpha|=\sqrt{\langle\alpha,\alpha\rangle}:$ the norm
of $ {R}^n,\forall \alpha\in {R}^n.$

$\bullet$~~$\langle A,B\rangle=tr(AB^T):$ the inner product in $
{R}^{n\times m},\forall A,B\in  {R}^{n\times m}.$

$\bullet$~~$|A|=\sqrt{tr(AA^T)}:$ the norm of $ {R}^{n\times
m},\forall A\in  {R}^{n\times m}$.

$\bullet$~~$l^2$: the space of all real-valued sequences
$x=(x_n)_{n\geq 0}$ satisfying
$$\|x\|_{l^2}\leq\sqrt{\displaystyle \sum_{i=1}^\infty x_i^2}<+\infty.$$

$\bullet$~~$l^2(H):$ the space of all H-valued sequence
$f=\{f^i\}_{i\geq 1}$ satisfying
$$\|f\|_{l^2(H)}\leq\sqrt{\displaystyle\sum_{i=1}^\infty||f^i||_H^2}<+\infty.$$

$\bullet$~~$l_{\mathscr{F}}^2(0, T, H):$ the space of all
$l^2(H)$-valued and ${\mathscr{F}}_t$-predictable processes
$f=\{f^i(t,\omega),\ (t,\omega)\in[0,T]\times\Omega\}_{i\geq1}$
satisfying
$$\|f\|_{l_{\mathscr{F}}^2(0, T, H)}\leq\sqrt{E\displaystyle\int_0^T\sum_{i=1}^\infty||f^i(t)||_H^2dt}<\infty.$$

$\bullet$~~$M_{\mathscr{F}}^2(0,T;H):$ the space of all $H$-valued
and ${\mathscr{F}}_t$-adapted processes $f=\{f(t,\omega),\
(t,\omega)\in[0,T]\times\Omega\}$ satisfying
$$\|f\|_{M_{\mathscr{F}}^2(0,T;H)}\leq\sqrt{E\displaystyle\int_0^T\|f(t)\|_H^2dt}<\infty.$$

$\bullet$~~$S_{\mathscr{F}}^2(0,T;H):$ the space of all $H$-valued
and ${\mathscr{F}}_t$-adapted  c\`{a}dl\`{a}g processes
$f=\{f(t,\omega),\ (t,\omega)\in[0,T]\times\Omega\}$ satisfying
$$\|f\|_{S_{\mathscr{F}}^2(0,T;H)}\leq\sqrt{E\displaystyle\sup_{0 \leq t \leq T}\|f(t)\|_H^2dt}<+\infty.$$

$\bullet$~~$L^2(\Omega,{\mathscr{F}},P;H):$ the space of all
$H$-valued random variables $\xi$ on $(\Omega,{\mathscr{F}},P)$
satisfying
$$\|\xi\|_{L^2(\Omega,{\mathscr{F}},P;H)}\leq E\|\xi\|_H^2<\infty.$$

The coefficients of the state equation (\ref{eq:2.1}) and the cost functional (\ref{eq:2.2}) are
defined as follows
 $$\xi: \Omega\longrightarrow {R}^n,$$
$$f:[0,T]\times \Omega\times
{R}^n\times {R}^{n\times d}\times
l^2({R}^n)\times U_1\times U_2\longrightarrow {R}^n,$$
$$l:[0,T]\times \Omega\times
{R}^n\times {R}^{n\times d}\times
l^2({R}^n)\times U_1\times U_2\longrightarrow {R}^1$$ and
$$\phi:\Omega\times {R}^{n}\longrightarrow {R}^1$$

Throughout this paper, we introduce the following basic assumptions
on coefficients $(\xi, f, l, \phi)$.
\begin{ass}\label{ass:2.1}   $\xi\in L^2(\Omega,{\mathscr{F}}_T,P;
 {R}^n)$ and the random mapping $f$ is
${\mathscr{P}}\bigotimes {\mathscr B}({R}^n)\bigotimes {\mathscr B}
({R}^{n\times d})\bigotimes {\mathscr B}(l^2({R}^n))\bigotimes
\\{\mathscr B}(U_1)\bigotimes
{\mathscr B}(U_2)$ measurable with $f(\cdot,0,0,0,0,0)\in M^2_{\mathscr
F}(0,T;{R}^n)$. For almost all $(t, \omega)\in [0, T]\times
\Omega$, $f(t,\omega, y,p,z,u_1,u_2)$ is Fr\'{e}chet differentiable w.r.t.
$(y,p,z,u_1,u_2)$ and the corresponding Fr\'{e}chet  derivatives $f_y,
f_p, f_z, f_{u_1},f_{u_2}$ are continuous and uniformly bounded.
\end{ass}
\begin{ass}\label{ass:2.2}
The random mapping $l$ is  ${\mathscr{P}}\bigotimes {\mathscr
B}({R}^n)\\\bigotimes {\mathscr B} ({R}^{n\times d})\bigotimes
{\mathscr B}(l^2({R}^n))\bigotimes {\mathscr B}(U_1)\bigotimes {\mathscr B}(U_2)$ measurable and
for almost all $(t, \omega)\in [0, T]\times \Omega$, $l$ is
Fr\'{e}chet differentiable w.r.t. $(y,p,z,u_1,u_2)$ with continuous
Fr\'{e}chet derivatives $l_y, l_q, l_z, l_{u_1},l_{u_2}$. The random mapping
$\phi$ is measurable and for almost all $(t,\omega)\in [0,
T]\times\Omega$, $\phi$ is Fr\'{e}chet differentiable w.r.t. $y$
with continuous Fr\'{e}chet derivative $\phi_y$. Moreover, for
almost all $(t,\omega)\in [0, T]\times \Omega$, there exists a
constant $C$ s.t. for all $(p,q,z,u_1,u_2)\in
 {R}^n\times {R}^{n\times d}\times
l^2( {R}^n)\times U_1 \times U_2$,
$$|l|\leq C(1+|y|^2+|q|^2+|z|^2+|u|_1^2+|u|_2^2),\ \ |\phi|\leq C(1+|y|^2),$$
$$|l_y|+|l_q|+|l_z|+|l_{u_1}|+|l_{u_2}|\leq C(1+|y|+|q|+|z|+|u_1|+|u_2|)$$
and $$ \ |\phi_y|\leq C(1+|y|).$$
\end{ass}
Under Assumption \ref{ass:2.1}, we can get from Lemma 2.3 in [7]
that for each $\big(u_1(\cdot), u_2(\cdot))\in {\cal A}_1\times{\cal A}_2$, the system (\ref{eq:2.1})
admits a unique strong solution. Furthermore, by Assumption \ref{ass:2.2} and a priori estimate
for BSDE driven by teugel martingales (see lemma 3.2 in [7]), it is
easy to check that
$$ |J(u_1(\cdot), u_2(\cdot))|<\infty.$$ So the Problem (P) is well-defined.

 \section{A PARTIAL INFORMATION SUFFICIENT MAXIMUM PRINCIPLE}
 \setcounter{equation}{0}
\renewcommand{\theequation}{3.\arabic{equation}}
In this section we want to study the  sufficient maximum principle
for Problem (P).

In our setting  the Hamiltonian function $H:[0,T]\times{R}^n\times
 {R}^{n\times d}\times l^2({R}^n)\times U_1\times U_2\times {R}^n \rightarrow {R}^1$  gets the following form:
\begin{equation}\label{eq:4.2}
\begin{array}{ll} \displaystyle
H(t,y,q,z,u_1,u_2,k)=&\langle k, -f(t,y,q,z,u_1,u_2)\rangle
\\&+l(t,y,q,z,u_1,u_2)
\end{array}
\end{equation}
The  adjoint  equation which fits into  the system (2.1) corresponding
   to the given admissible pair $({u}_1(\cdot), {u}_2(\cdot));  y(\cdot), q(\cdot), z(\cdot))$ is given by the following
   forward stochastic differential equation driven by multi-dimensional
Brownian motion $W$ and Teugel's martingales $\{H^i\}_{i=1}^\infty$:
   \begin{equation}\label{eq:3.2}
\left\{\begin{array}{ll}
dk(t)=-H_y(t,{y}(t),{q}(t),{z}(t),u_1(t),u_2(t),k(t))dt\\
\ \ \ \ \ \ \ \ \ \ \ -\displaystyle\sum_{i=1}^d H_{q^{i}}
(t,{y}(t),\bar{q}(t),{z}(t),u_1(t),u_2(t),k(t))dW^i(t)\\
\ \ \ \ \ \ \ \ \ \ \ -\displaystyle\sum_{i=1}^\infty H_{z^i}(t,{y}(t),{q}(t),{z}(t),u_1(t),u_2(t),k(t))dH^i(t)\\
k(0)=-\phi_y({y}(0)).
\end{array}\right.
\end{equation}
Under Assumptions \ref{ass:2.1}-\ref{ass:2.2},  the forward
stochastic differential equation (\ref{eq:3.2}) has a unique
solution $k(\cdot)\in {\cal S}^2_{\mathscr{F}}(0,T; {R}^n)$ by Lemma
2.1 in [7].

We now coming to a verification theorem for Problem (P).

{\bf Theorem 3.1 (A partial information sufficient maximum
principle)}  Let Assumptions \ref{ass:2.1}-\ref{ass:2.2} holds. Let $(\bar{u}_1(\cdot), \bar{u}_2(\cdot);
\bar{y}(\cdot),\bar{q}(\cdot),\bar{z}(\cdot))$ be an admissible pair
and $ \bar {k}(\cdot)$ be  the  unique strong solution of the
corresponding adjoint  equation (\ref{eq:3.2}).
Suppose   that the Hamiltonian function $H$ satisfies the following conditional
mini-maximum principle:
\begin {equation}   \label{eq:3.9}
\begin{array}{ll}
&\displaystyle\inf_{u_1\in{\cal U}_1}
E[H(t,\bar{y}(t),\bar{q}(t),\bar{z}(t),u_1,\bar{u}_2(t),\bar{k}(t))|{\cal G}_t]
\\&=E[H(t,\bar{y}(t),\bar{q}(t),\bar{z}(t),\bar{u}_1(t),\bar{u}_2(t),
\bar{k}(t))|{\cal G}_t]\\&=\displaystyle\sup_{u_2\in
{\cal
U}_2}E[H(t,\bar{y}(t),\bar{q}(t),\bar{z}(t),\bar{u}_1(t),u_2,\bar{k}(t))|{\cal G}_t].
\end{array}
\end {equation}
(i) Suppose that, for all $t\in [0,T],$  $\phi(y)$ is convex in $y$,  and
$$
(y,q,z,u_1)\mapsto
H(t,y,q,z,{u}_1,\bar{u}_2(t),\bar{k}(t))
$$
is convex.  Then,
$$
J(\bar{u}_1(\cdot),\bar{u}_2(\cdot))\leq
J(u_1(\cdot),\bar{u}_2(\cdot)),~~~for~~all ~~u_1(\cdot)\in{\cal
A}_1,
$$
and
$$
J(\bar{u}_1(\cdot),\bar{u}_2(\cdot))=\displaystyle\inf_{u_1(\cdot)\in{\cal
A}_1}J(u_1(\cdot),\bar{u}_2(\cdot)).
$$
(ii) Suppose that, for all $t\in [0,T],$ $\phi(y)$ is concave in $y$  and
$$
(y,q,z,u_2)\mapsto
H(t,y,q,z,\bar{u}_1(t),u_2,\bar{k}(t))
$$
is concave. Then,  $$ J(\bar{u}_1(\cdot),\bar{u}_2(\cdot))\geq
J(\bar{u}_1(\cdot),{u}_2(\cdot)),~~~for~~all ~~u_2(\cdot)\in{\cal
A}_2,
$$
and
$$
J(\bar{u}_1(\cdot),\bar{u}_2(\cdot))=\displaystyle\sup_{u_2(\cdot)\in{\cal
A}_2}J(\bar{u}_1(\cdot),u_2(\cdot)).
$$

(iii) If both case (i) and (ii) hold (which implies, in particular,
that $\phi(y)$  is an affine function), then
$(\bar{u}_1(\cdot),\bar{u}_2(\cdot))$ is an open-loop saddle point
and
\begin {equation}
\begin{array}{ll}
&~~~J(\bar{u}_1(\cdot),\bar{u}_2(\cdot))=\displaystyle\sup_{u_2(\cdot)\in{\cal
A}_2} (\inf_{u_1(\cdot)\in{\cal
A}_1}J(u_1(\cdot),u_2(\cdot)))\\
&=\inf_{u_1(\cdot)\in{\cal A}_1}(\sup_{u_2(\cdot)\in{\cal
A}_2}J(u_1(\cdot),u_2(\cdot))).
\end{array}
\end {equation}

{\bf Proof}:  (i) In the following,  we  consider a stochastic optimal control problem over ${\cal A}_1$ where the system  is
{\small \begin{equation} \label{eq:3.5}
\begin{array}{ll}
&y(t)=\xi+\displaystyle\int_t^Tf(s,y{(s)},q(s),z(s),u_1(t),\bar u_2(t))ds\\&
-\displaystyle\sum_{i=1}^d\int_t^Tq^i(s)dW^i(s)
-\displaystyle\sum_{i=1}^\infty\int_t^T z^i(s)dH^i(s),\ t\in[0,T]
\end{array}
\end {equation}
with the cost functional
\begin{equation}\label{eq:3.6}
\begin{array}{ll}
&J(u_1(\cdot),\bar u_2(\cdot))
\\&=E\displaystyle\bigg[\int_0^Tl(t,y(t),q(t),z(t),u_1(t),\bar u_2(t))dt
+\phi(y(0))\bigg].
\end{array}
\end{equation}}

 Our  optimal control problem is minimize
$J(u_1(\cdot), \bar{u}_2(\cdot))$ over $u_1(\cdot) \in {\cal A }_1$,
i.e., find  $\bar{u}_1(\cdot)\in {\cal A }_1$ such
that\begin{equation} J(\bar{u}_1(\cdot),
\bar{u}_2(\cdot))=\displaystyle\inf_{u_1(\cdot)\in {\cal
A}_1}J(u_1(\cdot), \bar{u}_2(\cdot)).
\end{equation}
Then for this case, it is easy to check that the Hamilton is  $H(t,y,q,z,u_1, \bar u_2(t),k)$
and  for the admissible control
$\bar u_1(\cdot)\in {\cal A}_1,$   the corresponding  sate process and the adjoint process
is still $(\bar{y}(t),\bar{q}(t),\bar{z}(t))$ and $\bar{k}(t),$  respectively.
And the optimality condition  is
\begin {equation}
  \begin{array}{ll}
  &\displaystyle\inf_{u_1\in{\cal U}_1}
  E[H(t,\bar{y}(t),\bar{q}(t),\bar{z}(t),u_1,\bar{u}_2(t),\bar{k}(t))|{\cal G}_t]
  \\&=E[H(t,\bar{y}(t),\bar{q}(t),\bar{z}(t),\bar{u}_1(t),\bar{u}_2(t)|{\cal G}_t].
  \end{array}
  \end {equation}

 Thus from the partial information sufficient maximum principle for optimal control (see Theorem
 3.1 in [9]), we conclude that $\bar{u}_1(\cdot)$ is  the optimal
 control of the optimal control problem, i.e.,
$$
J(\bar{u}_1(\cdot),\bar{u}_2(\cdot))\leq
J(u_1(\cdot),\bar{u}_2(\cdot)),~~~for~~all ~~u_1(\cdot)\in{\cal
A}_1,
$$
and
$$
J(\bar{u}_1(\cdot),\bar{u}_2(\cdot))=\displaystyle\inf_{u_1(\cdot)\in{\cal
A}_1}J(u_1(\cdot),\bar{u}_2(\cdot)).
$$
     The proof of (i)is complete  .

(ii)This statement can be proved in a similar way as (i).\\
(iii) if both (i) and (ii) hold, then
$$
J(\bar{u}_1(\cdot),u_2(\cdot))\leq
J(\bar{u}_1(\cdot),\bar{u}_2(\cdot))\leq
J(u_1(\cdot),\bar{u}_2(\cdot)),
$$
for any $(u_1(\cdot),u_2(\cdot))\in {\cal A}_1\times {\cal A}_2.$
Thereby,
$$
\begin{array}{ll}
J(\bar{u}_1(\cdot),\bar{u}_2(\cdot))&\leq
\displaystyle\inf_{u_1(\cdot)\in {\cal
A}_1}J(u_1(\cdot),\bar{u}_2(\cdot))\\
&\leq \sup_{u_2(\cdot)\in {\cal
A}_2}(\displaystyle\inf_{u_1(\cdot)\in {\cal
A}_1}J(u_1(\cdot),u_2(\cdot))).
\end{array}
$$
On the other hand,
$$
\begin{array}{ll}
J(\bar{u}_1(\cdot),\bar{u}_2(\cdot))&\geq
\displaystyle\sup_{u_2(\cdot)\in {\cal
A}_2}J(\bar{u}_1(\cdot),{u}_2(\cdot))\\
&\geq \displaystyle\inf_{u_1(\cdot)\in {\cal
A}_1}(\displaystyle\sup_{u_2(\cdot)\in {\cal
A}_2}J(u_1(\cdot),u_2(\cdot))).
\end{array}
$$
Now, due to the inequality
$$
\begin{array}{ll}
&\displaystyle\inf_{u_1(\cdot)\in {\cal A}_1}(\sup_{u_2(\cdot)\in
{\cal A}_2}J(u_1(\cdot),u_2(\cdot)))\\
\geq &\sup_{u_2(\cdot)\in {\cal A}_2}(\inf_{u_1(\cdot)\in {\cal
A}_1}J(u_1(\cdot),u_2(\cdot))),
\end{array}
$$
we have
$$
\begin{array}{ll}
J(\bar{u}_1(\cdot),\bar{u}_2(\cdot))&=\displaystyle\sup_{u_2\in\Pi}(\inf_{u_1\in\Theta}J(u_1,u_2))\\
&=\inf_{u_1\in\Theta}(\sup_{u_1\in\Pi}J(u_1,u_2)).
\end{array}
$$
The proof of Theorem 3.1 is completed.

If the control process $(u_1(\cdot),u_2(\cdot))$ is admissible
adopted to the filtration ${\cal F}_t$ we have the following
full information sufficient maximum principle.

{\bf Corollary 3.1 } Suppose that ${\cal G}_t={\cal F}_t$  . Moreover, suppose that, for all $t\in [0,T],$
the following maximum principle holds:
\begin {equation}
\begin{array}{ll}
&\displaystyle\inf_{u_1\in{\cal U}_1}
H(t,\bar{x}(t),\bar{y}(t),\bar{z}(t),u_1,\bar{u}_2(t),\bar{p}(t),\bar{q}(t),\bar{k}(t))\\
=&H(t,\bar{x}(t),\bar{y}(t),\bar{z}(t),\bar{u}_1(t),\bar{u}_2(t),\bar{p}(t),\bar{q}(t),\bar{k}(t))
\\
=&\displaystyle\sup_{u_2\in {\cal
U}_2}H(t,\bar{x}(t),\bar{y}(t),\bar{z}(t),\bar{u}_1(t),u_2,\bar{p}(t),\bar{q}(t),\bar{k}(t)).
\end{array}
\end {equation}

(i) Suppose that, for all $t\in [0,T],$ $\phi(y)$ is convex in $y$  and
$$
(x,y,z,u_1)\mapsto
H(t,x,y,z,{u}_1,\bar{u}_2(t),\bar{p}(t),\bar{q}(t),\bar{k}(t))
$$
is convex.  Then,
$$
J(\bar{u}_1(\cdot),\bar{u}_2(\cdot))\leq
J(u_1(\cdot),\bar{u}_2(\cdot)),~~~for~~all ~~u_1(\cdot)\in{\cal
A}_1,
$$
and
$$
J(\bar{u}_1(\cdot),\bar{u}_2(\cdot))=\displaystyle\inf_{u_1(\cdot)\in{\cal
A}_1}J(u_1(\cdot),\bar{u}_2(\cdot)).
$$
(ii) Suppose that, for all $t\in [0,T],$  $\phi(y)$ is concave in $y,$  and
$$
(x,y,z,u_2)\mapsto
H(t,x,y,z,\bar{u}_1(t),u_2,\bar{p}(t),\bar{q}(t),\bar{k}(t))
$$
is concave. Then,  $$ J(\bar{u}_1(\cdot),\bar{u}_2(\cdot))\geq
J(\bar{u}_1(\cdot),{u}_2(\cdot)),~~~for~~all ~~u_2(\cdot)\in{\cal
A}_2,
$$
$$
J(\bar{u}_1(\cdot),\bar{u}_2(\cdot))=\displaystyle\sup_{u_2(\cdot)\in{\cal
A}_2}J(\bar{u}_1(\cdot),u_2(\cdot)).
$$

(iii) If both case (i) and (ii) hold (which implies, in particular,
that $\phi(y)$ is an affine function), then
$(\bar{u}_1(\cdot),\bar{u}_2(\cdot))$ is an open-loop saddle point
based on the information flow ${\cal F}_t$ and
\begin {equation}
\begin{array}{ll}
J(\bar{u}_1(\cdot),\bar{u}_2(\cdot))&=\displaystyle\sup_{u_2(\cdot)\in{\cal
A}_2} (\inf_{u_1(\cdot)\in{\cal
A}_1}J(u_1(\cdot),u_2(\cdot)))\\
&=\displaystyle\inf_{u_1(\cdot)\in{\cal
A}_1}(\sup_{u_2(\cdot)\in{\cal A}_2}J(u_1(\cdot),u_2(\cdot))).
\end{array}
\end {equation}


  \section{A PARTIAL INFORMATION NECESSARY MAXIMUM PRINCIPLE}
  \setcounter{equation}{0}
 \renewcommand{\theequation}{4.\arabic{equation}}

 In this section, we give a necessary maximum principle for Problem (P).

{\bf Theorem 4.1 (A partial information  necessary maximum
principle)}
 Under Assumptions \ref{ass:2.1}-\ref{ass:2.2}, let $(\bar{u}_1(\cdot), \bar{u}_2(\cdot))$
  be an optimal pair of  Problem (P).
Suppose that  $(\bar{y}(\cdot), \bar{q}(\cdot), \bar{z}(\cdot))$ is
 the state process of the system (\ref{eq:2.1}) corresponding to
 the admissible control $(\bar{u}_1(\cdot), \bar{u}_2(\cdot)).$
 Let $(\bar k(\cdot) )$ be  the unique
 solution of the  adjoint equation (\ref{eq:3.2}) corresponding $
 (\bar{u}_1(\cdot), \bar{u}_2(\cdot);  \bar{y}(\cdot), \bar{q}(\cdot), \bar{z}(\cdot))$.
 Then for $i=1,2,$ we have
\begin{equation}\label{eq:4.1}
\begin{array}{ll}
\big\langle E[H_{u_i}(t)|{\cal G}_t], u_i-\bar{u}_i(t)\big\rangle \geq 0, \forall
 u_i \in U_i a.s.a.e.,
\end{array}
\end{equation}

where
$$H_{u_i}(t)=:H_{u_i}(t,\bar{y}(t),\bar{q}(t),\bar{z}(t),\bar u_1(t),\bar{u}_2(t),\bar{k}(t)).$$
 {\bf Proof:} Since  $(\bar{u}_1(\cdot),\bar{u}_2(\cdot))$ is an
optimal open-loop control, then  $(\bar{u}_1(\cdot),\bar{u}_2(\cdot))$
is an open-loop
saddle point, i.e.,
\begin{eqnarray}
J(\bar{u}_1(\cdot), {u}_2(\cdot))\leq J(\bar{u}_1(\cdot),
\bar{u}_2(\cdot))\leq J(u_1(\cdot), \bar{u}_2(\cdot)).
\end{eqnarray}
So we have \begin{equation}\label{eq:3.3}
    J_1(\bar{u}_1(\cdot),
\bar{u}_2(\cdot))=\displaystyle\min_{u_1(\cdot)\in {\cal
A}_1}J_1(u_1(\cdot), \bar{u}_2(\cdot)), \end{equation}
and
\begin{equation}\label{eq:3.4}
 J_2(\bar{u}_1(\cdot),
\bar{u}_2(\cdot))=\displaystyle\max_{u_2(\cdot)\in {\cal
A}_2}J_2(\bar u_1(\cdot), {u}_2(\cdot)).
\end{equation}
By (\ref{eq:3.3}),  $\bar u_1(\cdot)$  can be regarded as
an optimal control  of the
optimal control problem where
the controlled system is
(\ref{eq:3.5}) and the cost functional
is (\ref{eq:3.6}).
Then for this case, it is easy to check that the Hamilton is  $H(t,y,q,z,u_1, \bar u_2(t),k)$
and  for the optimal control
$\bar u_1(\cdot)\in {\cal A}_1,$   the corresponding optimal  sate process and the adjoint process
is still $(\bar{y}(t),\bar{q}(t),\bar{z}(t))$ and $\bar{k}(t),$  respectively.
Thus applying the     partial
necessary stochastic maximum principle for  optimal control problems
(see Theorem 4.1 in [9]), we can obtain \ref{eq:4.1} for $i=1$.
Similarly, from (\ref{eq:3.4}), we can obtain (\ref{eq:4.1}) for $i=2$.

 The proof is complete.

\section {An EXAMPLE:  A Linear quadratic problem}

In this section, we will apply our stochastic maximum principles to
a Linear quadratic problem under partial information, i.e. minimizing the
following quadratic cost functional over $u(\cdot) \in \cal A$:
\begin{equation}\label{eq:6.1}
\begin{array}{ll}
J(u_1(\cdot),u_2(\cdot)):=&E\langle M,y(0)\rangle
+E\displaystyle\int_0^T\langle E(s),y(s)\rangle ds
\\&+\sum_{i=1}^dE\displaystyle\int_0^T\langle
F^i(s),q^i(s)\rangle ds\\&+ \sum_{i=1}^\infty
E\displaystyle\int_0^T\langle G^i(s),z^i(s)\rangle ds
\\&+E\displaystyle\int_0^T\langle N_1(s)u_1(s),u_1(s)\rangle ds
\\&-E\displaystyle\int_0^T\langle N_2(s)u_2(s),u_2(s)\rangle ds,
\end{array}
\end{equation}
where the state process $(y(\cdot), q(\cdot), z(\cdot))$ is the
solution to the controlled linear backward stochastic system below:
\begin{equation}\label{eq:6.2}
\left\{\begin{array}{lll}{}
dy(t)=-\bigg[A(t)y(t)+\displaystyle\sum_{i=1}^dB^i(t)q^i(t)
\\+\displaystyle\sum_{i=1}^\infty C^{i}(t)z^i(t)+D_1(t)u_1(t)+D_2(t)u_2(t)\bigg]dt\\
 +\displaystyle\sum_{i=1}^dq^i(t)dW^i(t)+\displaystyle\sum_{i=1}^\infty z^i(t)dH^{i}(t),\\
y(T)=\xi.
\end{array}
\right.
\end {equation}
This problem we denote by Problem (LQ).
To study this problem, we need the assumptions on the coefficients
as follows.
\begin{ass}\label{ass:5.1}
The matrix-valued functions  $A:  [0,T]  \rightarrow {R}^{n\times n};
B^i: [0,T] \rightarrow {R}^{n\times n},i=1,2,\cdots,d; C^i:  [0,T]
\rightarrow {R}^{n\times n},i=1,2,\cdots;  D^i:  [0,T] \rightarrow
{R}^{n\times m_i},i=1,2; E:  [0,T] \rightarrow {R}^{n\times n}; F^i:  [0,T]
\rightarrow       {R}^{n\times m};
 {R}^{n}, i=1,2,\cdots d, G^i:  [0,T]
\rightarrow  {R}^{n}, i=1,2,\cdots; N^i:[0,T] \rightarrow
{R}^{m_i\times m_i},i=1,2 $
 and the  matrix $M\in {R}^{n}$ are
uniformly bounded.          Moreover, $N^i$ is  uniformly positive,
i.e. $ N^i\geq \delta I, (i=1,2) $ for   some positive constant   $\delta$  a.s..
\end{ass}

\begin{ass}\label{ass:5.3}
There is no further constraint imposed on the control processes, the set all admissible control processes is
$$
  \begin{array}{ll}
&{\cal A}_1\times {\cal A}_2=\bigg\{(u_1(\cdot),u_2(\cdot)):(u_1(\cdot),u_2(\cdot))\ is\ {\cal G}_t
\\&~~~~~~~-predictable\ process \
with\ values\ in\  {R}^{m_1}\times {R}^{m_2}\\&~~~~~~ \ and\
E\displaystyle\int_0^T|u(t)|^2dt< \infty \bigg\}.
  \end{array}
$$

\end{ass}

In what follows, we will utilize the stochastic maximum principle to
study the dual representation of the optimal control to Problem (LQ).
%

We first define the Hamiltonian function $H:\Omega\times[0,T] \times
{R}^n\times  {R}^{n\times d}\times l^2(  {R}^n)\times U_1\times U_2\times
{R}^n\rightarrow   {R}^1$ by
\begin{eqnarray}\label{eq:6.7}
&& H(t,y,q,z,u_1,u_2,k)\nonumber\\&&=-\bigg\langle k,
A(t)y+\displaystyle\sum_{i=1}^dB^i(t)q^i+\displaystyle\sum_{i=1}^\infty
C^{i}(t)z^i\nonumber \\&&+D_1(t)u_1+D_2(t)u_2\bigg\rangle\nonumber +\langle E(t), y\rangle \\&&\nonumber +\displaystyle\sum_{i=1}^d\langle
F^i(t),q^i\rangle\nonumber+\langle N_1(t)u_1, u_1\rangle\\&&~~~~
+\displaystyle\sum_{i=1}^\infty \langle G^i(t),z^i\rangle
  -  \langle N_2(t)u_1, u_1\rangle.
\end{eqnarray}
Then the adjoint equation
corresponding to an admissible pair $(u_1(\cdot), u_2(\cdot); y(\cdot),q(\cdot),z(\cdot)))$
can be rewritten as
\begin{equation}\label{5.8}
\displaystyle\left\{\begin{array}{lll}{} dk(t)=-H_y(t, {y}(t),
{q}(t), {z}(t), {u}_1(t), {u}_2(t),k(t))dt\\ \ \ \ \ \
-\displaystyle\sum_{i=1}^d H_q^{i} (t, {y}(t), {q}(t),
{z}(t), {u}_1(t), {u}_2(t),k(t))dW^i(t)\\
\ \ \ \ \ \  -\displaystyle\sum_{i=1}^\infty
H_{z^i}(t, {y}(t-), {q}(t), {z}(t), {u_1}(t), {u_2}(t),k(t))dH^i(t),\\
k(0)=-M.
\end{array}
\right.
\end{equation}

Under Assumption \ref{ass:5.1}, for any admissible pair
$({u_1}(\cdot),{u_2}(\cdot); {y}(\cdot), {q}(\cdot), {z}(\cdot))$, the adjoint
equation (\ref{5.8}) has a unique solution $k(\cdot)$ in view of
Lemma 2.1 in [9].

It is time to give the  the dual characterization of the optimal
control.

   \begin{thm}
Let  Assumptions \ref{ass:5.1}-\ref{ass:5.3} be satisfied.
 Then, a necessary and
sufficient condition for an admissible pair $(u_1(\cdot),u_2(\cdot); y(\cdot),
q(\cdot), z(\cdot))$ to be an optimal pair of  Problem (LQ) is the
control $(u_1(\cdot), u_2(\cdot))$ has
  the representation
\begin{equation}\label{eq:3.20}
\begin{array}{ll}
    u_1(t) =&-\frac{1}{2}N^{-1}_1(t)D_1^*(t)E[k(t)|{\cal G}_t],  \\
       u_2(t) =&\frac{1}{2}N^{-1}_2(t)D_2^*(t)E[k(t)|{\cal G}_t],
     \end{array}
\end{equation}
where $k(\cdot)$ is the unique
solution of the adjoint equation (\ref{5.8}) corresponding to the admissible pair
$({u_1}(\cdot),u_2(\cdot); {y}(\cdot), q(\cdot), {z}(\cdot))$.
     \end{thm}

{\bf Proof:}For the necessary part, let $(u_1(\cdot),u_2(\cdot); y(\cdot), q(\cdot), \\z(\cdot))$ be an
optimal pair,  then by the
necessary optimality condition  (\ref{eq:4.1}) and $U_{i}=R^{m_i}(i=1,2)$, we have
$$H_{u_i}(t,{y}(t), q(t),{z}(t),  {u_1}(t), u_2(t), k(t))=0,\ \ a.e.\ a.s..$$
Noticing the definition of $H$ in (\ref{eq:6.7}), we get
$$2N_1(t)u_1(t)+D_1^*(t)E[k(t)|{\cal G}_t]=0,\ \ a.e.\ a.s.$$
and       $$-2N_2(t)u_2(t)+D^*_2(t)E[k(t)|{\cal G}_t]=0,\ \ a.e.\ a.s.$$
So  the optimal control $(u_1(\cdot),u_2(\cdot))$ has the dual presentation
(\ref{eq:3.20}).\\
For the sufficient part, let $(u_1(\cdot),u_2(\cdot); y(\cdot), q(\cdot),
z(\cdot))$ be an admissible pair satisfying (\ref{eq:3.20}). By the
classical technique of completing squares,  from (\ref{eq:3.20}), we
can claim that $(u_1(\cdot),u_2(\cdot); y(\cdot), q(\cdot),z(\cdot))$ satisfies the
optimality condition (\ref{eq:3.9}) in Theorem 3.1. Moreover,
from Assumption \ref{ass:5.1}-\ref{ass:5.3}, it is easy
to check that  all other conditions in Theorem 3.1 are
satisfied. Hence $(u_1(\cdot),u_2(\cdot); y(\cdot), q(\cdot),z(\cdot))$ is
an optimal pair by Theorem 3.1.

\section{Conclusion}
In this paper we have proved partial information necessary
as well as sufficient conditions of  optimality for  the stochastic    differential
games  driven by  backward stochastic control systems associated with L\'{e}vy process.
     As an application, some  linear quadratic stochastic differential games  problem
      of backward stochastic system is discussed and the optimal control process is characterized
       explicitly by adjoint process.
       Our main results could be seen as an extension of the stochastic optimal control problem
     studied   in [9], to the stochastic differential games problem.

\end{document}